\documentclass{article}
\usepackage{amsmath,amssymb,amsthm}
\newtheorem{Thm}{Theorem}
\newtheorem{Lem}{Lemma}
\newtheorem{Prop}{Proposition}

\theoremstyle{definition}

\begin{document}

\title{Classification of double power nonlinear functions}
\author{Shinji Kawano}
\date{}
\maketitle

\begin{abstract}

On investigating the conditions of existence and uniqueness of positive solutions to the problem,
\begin{equation}
 \begin{cases}
   \triangle u+f(u)=0  & \text{in  $\mathbb{R}^n$},\\
   \displaystyle \lim_{\lvert x \rvert \to \infty} u(x)  =0, \label{Intro2}
 \end{cases}  
\end{equation}
examination of the function $f$ is necessary.

In fact for a double power nonlinear function
\begin{equation*}
f(u)=-\omega u + u^p - u^q, \qquad  \omega>0, ~~q>p>1,
\end{equation*}
a positive solution to \eqref{Intro2} exists if and only if $F(u):=\int_0^u f(s)ds >0$ for some $u>0$.
Moreover the solution is unique if $f$ satifies furthermore the condition
that $\tilde{f}(u):= (uf'(u))'f(u)-uf'(u)^2 <0$ for any $u>0$.

In the former papers we considered the relations between $f$, $F$, $\tilde{f}$ and $\tilde{F}$. 
In the present paper we remark that these relations also hold for more general nonlinear functions.
\end{abstract}

\section{Introduction}

We shall consider a boundary value problem
\begin{equation}
 \begin{cases}
   u_{rr}+ \dfrac{n-1}{r}u_r+f(u)=0 & \text{for $r>0$}, \\
   u_r(0)=0, \\
   \displaystyle \lim_{r \to \infty} u(r) =0, \label{b}
 \end{cases} 
\end{equation} 
where $n \in \mathbb{N}$ and
\begin{equation*}
f(u)=-\omega u + u^p - u^q, \qquad \omega>0, ~~q>p>1. 
\end{equation*}              
The above problem arises in the study of 
\begin{equation}
 \begin{cases}
   \triangle u+f(u) =0  & \text{in  $\mathbb{R}^n$},\\
   \displaystyle \lim_{\lvert x \rvert \to \infty} u(x)  =0. \label{a}
 \end{cases}  
\end{equation}
Indeed, the classical work of Gidas, Ni and Nirenberg~\cite{G1,G2} tells us that any positive solution to \eqref{a} is radially symmetric. 
On the other hand, for a solution $u(r)$ of \eqref{b}, $v(x):=u(\lvert x \rvert)$ is a solution to \eqref{a}.  

The condition to assure the existence of positive solutions to \eqref{a} (and so \eqref{b}) was given by 
Berestycki and Lions~\cite{B1} and Berestycki, Lions and Peletier~\cite{B2}:

\begin{Prop}
A positive solution to \eqref{b} exists if and only if 
\begin{equation}
F(u):=\int_0^u f(s)ds >0, \qquad \text{for some} \quad u>0.       \label{existence}
\end{equation} 
\end{Prop}     

Uniqueness of positive solutions to \eqref{b} had long remained unknown.
Finally in 1998 Ouyang and Shi~\cite{OS} proved uniqueness for \eqref{b} with $f$ satisfying the additional condition
(See also Kwong and Zhang~\cite{KZ}):

\begin{Prop}
If $f$ satisfies furthermore the following condition, then the positive solution is unique;
\begin{equation}
\tilde{f}(u):= (uf'(u))'f(u)-uf'(u)^2 <0, \qquad \text{for any} \quad u>0.      \label{unique}
\end{equation}     
\end{Prop}

In the former paper~\cite{Kawano3} we have presented the following results:
\begin{Thm}
The existence condition~\eqref{existence} is equivalent to the following condition;
\begin{equation}
\tilde{F}(u)= (uf(u))'F(u)-uf(u)^2 <0, \qquad \text{for any} \quad u>0.      \label{exitilde}
\end{equation} 

The uniqueness condition~\eqref{unique} is equivalent to the following condition;
\begin{equation}
f(u) >0, \qquad \text{for some} \quad u>0.       \label{uniquetilde}
\end{equation} 
\end{Thm}   \label{thm}

In the present paper we remark that these relations also hold for more general nonlinear functions.

This paper is organized as follows. In section 2, we give the statement of the main result.
In section~3, we give a simpler proof of Theorem~\ref{thm}, as an application of the formula given in the main result.

\section{Main result}

To state the main result, from now on we mean a more general function
\begin{equation}
 f(u)=-au^p+bu^q-cu^r,\qquad \text{for}  ~u>0,    \label{general}
\end{equation} 
where $a,b,c>0$ and $p<q<r$, by the same notation~$f$.
We also steal the former notation 
\begin{equation}
\tilde{f}(u)= (uf'(u))'f(u)-uf'(u)^2 .      \label{tilde}
\end{equation}    
Following is our main result.
\begin{Thm}
There only can occur the following three cases;
\begin{itemize}
\item[(a)]~~~~~ $a<b\dfrac{r-q}{r-p}\left[ \dfrac{b(q-p)}{c(r-p)} \right]^{\frac{q-p}{r-q}}$ $\Longleftrightarrow$ $f$ has positive parts  $\Longleftrightarrow$ $\tilde{f}$ remains negative.
\item[(b)]~~~~~ $a=b\dfrac{r-q}{r-p}\left[ \dfrac{b(q-p)}{c(r-p)} \right]^{\frac{q-p}{r-q}}$ $\Longleftrightarrow$ $f$ has just one zero $\Longleftrightarrow$ $\tilde{f}$ has just one zero.
\item[(c)]~~~~~ $a>b\dfrac{r-q}{r-p}\left[ \dfrac{b(q-p)}{c(r-p)} \right]^{\frac{q-p}{r-q}}$ $\Longleftrightarrow$ $f$ remains negative $\Longleftrightarrow$ $\tilde{f}$ has positive parts.
\end{itemize}
\end{Thm}    \label{th}

\begin{proof}
The statement with respect to $f$ is trivial.
We obtain from the definition~\eqref{tilde} that
\begin{equation*}
\tilde{f} = -ab(q-p)^2 u^{q+p-1} + ca(r-p)^2 u^{r+p-1} - bc(r-q)^2 u^{r+q-1} .
\end{equation*}
This is in the form of~\eqref{general} and use the result with respect to $f$.  
\end{proof}

\section{Application of the formula}

In this section we use the formula given in Theorem~\ref{th} to give a simpler proof of Theorem~\ref{thm}.  

The following two lemmas asserts the Theorem~\ref{thm}.
\begin{Lem}
Both the existence condition~\eqref{existence} and the condition~\eqref{exitilde} are equivalent to 
\begin{equation*}
\omega < \omega_{p,q},
\end{equation*}
 where 
\begin{equation*}
\omega_{p,q}=\dfrac{2(q-p)}{(p+1)(q-1)} \left[ \dfrac{(p-1)(q+1)}{(p+1)(q-1)} \right] ^{\frac{p-1}{q-p}}. 
\end{equation*}
(See Ouyang and Shi~\cite{OS} and the appendix of Fukuizumi~\cite{Fukuizumi}.) 
\end{Lem}

\begin{proof}
\begin{equation*}
F(u)= -\frac{\omega}{2}u^2+\frac{u^{p+1}}{p+1}-\frac{u^{q+1}}{q+1}, \qquad  \omega>0, ~~q>p>1,
\end{equation*}
is in the form of~\eqref{general} and use the result (a) of Theorem~\ref{th}. 
\end{proof}

\begin{Lem}
Both the uniqueness condition~\eqref{unique} and the condition~\eqref{uniquetilde} are equivalent to
\begin{equation*}
\omega < \eta_{p,q},
\end{equation*}
where
\begin{equation*}
\eta_{p,q}=\dfrac{q-p}{q-1}\left[ \dfrac{p-1}{q-1}\right]^{\frac{p-1}{q-p}}.
\end{equation*}
\end{Lem}

\begin{proof}
\begin{equation*}
f(u)=-\omega u + u^p - u^q, \qquad  \omega>0, ~~q>p>1,
\end{equation*}
is in the form of~\eqref{general} and use the result (a) of Theorem~\ref{th}. 
\end{proof}

\end{document}